\documentclass{amsart}

\usepackage{amsfonts}
\usepackage[leqno]{amsmath}
\usepackage[noadjust]{cite}
\usepackage{amssymb}
\usepackage{amstext}
\usepackage{caption}
\usepackage{subcaption}
\usepackage{graphicx}
\DeclareGraphicsExtensions{eps,jpg,tif}
\usepackage{titletoc}
\usepackage[french,english]{babel}
\usepackage{lineno}
\usepackage{calc}
\usepackage{mathabx}
\usepackage{pgf,tikz}
\usepackage[colorlinks=true,linkcolor=blue,pdftex=true,urlcolor=blue,hyperindex=true]{hyperref}	
\usepackage[nameinlink]{cleveref}

\newtheorem{theorem}{Theorem}[section]

\newtheorem{corollary}[theorem]{Corollary}

\newtheorem{example}[theorem]{Example}

\newtheorem{lemma}[theorem]{Lemma}
\newtheorem{notation}[theorem]{Notation}

\newtheorem{proposition}[theorem]{Proposition}
\newtheorem{remark}[theorem]{Remark}

\numberwithin{equation}{section}

\newcommand{\m}[1]{\left\vert #1\right\vert}

\newcommand{\C}{\mathbb C}
\newcommand{\R}{\mathbb R}
\newcommand{\Z}{\mathbb Z}
\newcommand{\Na}{\mathbb N}
\newcommand{\dsum}[3]{\displaystyle\sum_{#1}^{#2}#3}

\newcommand{\HolC}{\mathcal O(\C)}
\newcommand{\HolCstar}{\mathcal O(\C^*)}
\newcommand{\MerC}{\mathcal M(\C)}
\newcommand{\MerCstar}{\mathcal M(\C^*)}

\begin{document}

\title[Meromorphic solutions of linear $q$-difference equations]{Meromorphic solutions of linear $q$-difference equations}
\author{Alberto Lastra}
\author{Pascal Remy}
\address{University of Alcal\'{a}\\
Departamento de Fis\'{i}sica y Matem\'{a}ticas\\
Ap. de Correos 20, E-28871 Alcal\'{a} de Henares (Madrid), Spain}
\email{alberto.lastra@uah.es}
\address{Laboratoire de Math\'{e}matiques de Versailles\\
Universit\'{e} de Versailles Saint-Quentin\\
45 avenue des Etats-Unis, 78035 Versailles cedex, France}
\email{pascal.remy@uvsq.fr ; pascal.remy.maths@gmail.com}



\begin{abstract}
In this article, we construct explicit meromorphic solutions of first order linear $q$-difference equations in the complex domain and we describe the location of all their zeros and poles. The homogeneous case leans on the study of four fundamental equations, providing the previous informations in the framework of entire or meromorphic coefficients. The inhomogeneous situation, which stems from the homogeneous one and two fundamental equations, is also described in detail. We also address the case of higher-order linear $q$-difference equations, using a classical factorization argument. All these results are illustrated by several examples.
\end{abstract}

\maketitle


\section{Introduction}

The study of $q$-difference equations in the complex domain has experienced great interest in recent decades not only due to its inherent interest and numerous applications, but also motivated by and leaning on different previous theories.

On the one hand, one can point out different works based on Nevanlinna value distribution theory treating meromorphic solutions to $q$-difference equations. In~\cite{zhyi} (see also its references) certain properties of the image of $q$-difference operators acting meromorphic functions guarantee injectivity of the operator. In~\cite{zhtu}, estimates on the growth of the solutions to nonlinear $q$-difference equations are provided. Other recent advances based on Nevanlinna theory are~\cite{dugazhzh,wa,qiya}.

On the other hand, the study of solutions to (systems of) ordinary differential equations in the complex domain has also been considered in the $q$-analog framework. We only give some examples of trends and topics under study in this direction. This is the case of $q$-Gevrey asymptotic expansions and summability (see~\cite{rasazh,razh,taya,zh} among many others), Newton polygon techniques (see the recent work~\cite{cafo} and the references therein), summability techniques and tools such as $q$-analogs of Laplace transform to be applied in a Borel-Laplace methodology for providing analytic solutions to $q$-difference equations from formal ones~\cite{ta}, or the asymptotic study of $q$-difference-differential equations (see~\cite{lamasa,ma,ta2,ya} among others, and the references therein) have also been achieved. We also refer to the work~\cite{dr} where the author describes a procedure to provide meromorphic solutions to linear $q$-difference equations with rational coefficients of any positive order by means of $q$-analogs of Borel and Laplace transformations. Note that in~\cite{pr}, the existence of meromorphic solutions of general linear $q$-difference equations with meromorphic coefficients in $\C^*$ is proved via theoretical tools.

In the present work, we focus on the latter and we propose to build explicit solutions for them.

First, we consider the case of first-order equations of the form
\begin{equation}\label{EqIntro}
y(qx)=m(x)y(x)+r(x),
\end{equation}
where $m(x)$ and $r(x)$ are two meromorphic functions on $\C^*$ (Section \ref{sec21}). We start by determining explicit meromorphic solutions of six fundamental equations (Eqs. (\ref{BasicEquation1})-(\ref{BasicEquation6}), Section \ref{sec211}). This allows, by means of the Weierstrass Factorization Theorem and the Birkhoff Decomposition, to provide explicit meromorphic solutions of Eq. (\ref{EqIntro}) in the homogeneous case, that is when $r(x)\equiv0$ (Theorem \ref{HomogeneousHolC} and Corollaries \ref{HomogeneousMerC} and \ref{HomogeneousMerCstar}, Section \ref{sec212}). In particular, we prove that the location of the zeros and poles of such solutions is depending on their distance at the origin. The inhomogeneous case is treated in Section \ref{SectionInhomogeneousLinearQDiffOrder1} by means of a convenient change of variable and a Laurent-type series decomposition.

In Section \ref{ordern}, we conclude with some words on the consideration of linear $q$-difference equations of higher order, where the factorization known in the literature (see for instance \cite{Sau04,mazh,BG41}) allows to attain positive results. A complete study of such equations is left for a future research.

Some illustrative examples are also provided.

\begin{notation}
\textnormal{All along this work, we consider a nonzero complex number $q\in\C^*$ with $\m{q}>1$, and we use the following notations:
\begin{itemize}
\item $\Na$ (resp. $\Na^*$) the set of all the nonnegative (resp. positive) integers;
\item$\Z$ (resp. $\Z^*$) the set of all the integers (resp. nonzero integers);
\item$\C$ (resp. $\C^*$) the set of all the complex (resp. nonzero complex) numbers;
\item $\HolC$ the set of all the entire functions;
\item $\MerC$ (resp. $\MerCstar$) the set of all the meromorphic functions on $\C$ (resp. $\C^*$);
\item $q^{\Z}$ the discrete $q$-spiral $\{q^n;n\in\Z\}$;
\item $q^{\Na}$ (resp. $q^{\Na^*}$, $q^{-\Na}$, $q^{-\Na^*}$) the discrete $q$-half-spiral $\{q^n;n\in\Na\}$ (resp. $\Na^*$, $-\Na$, $-\Na^*$);
\item $\mathcal Z_m$ (resp. $\mathcal Z_m^*$) the set of all the zeros (resp. nonzero zeros) of a meromorphic function $m$;
\item$\mathcal Z^*_{m,\leq\rho}$ (resp. $\mathcal Z^*_{m,>\rho}$) with $\rho>0$ the set, possibly empty, of all the zeros $a\in\mathcal Z^*_m$ such that $\m{a}\leq\rho$ (resp. $\m{a}>\rho$);
\item $\mathcal P_m$ (resp. $\mathcal P_m^*$) the set of all the poles (resp. nonzero poles) of a meromorphic function $m$;
\item$\mathcal P^*_{m,\leq\rho}$ (resp. $\mathcal P^*_{m,>\rho}$) with $\rho>0$ the set, possibly empty, of all the poles $a\in\mathcal P^*_m$ such that $\m{a}\leq\rho$ (resp. $\m{a}>\rho$);
\item $P_{m,a}(x)$, with $a\in\mathcal P_m$, the principal part of $m$ at $a$, that is the polynomial in $1/(x-a)$ without constant term such that $m(x)-P_{m,a}(x)$ has a removable singularity at $a$;
\item$\mu_{m,a}$ the order of multiplicity of $a\in\mathcal Z_m^*\cup\mathcal P_m^*$.
\end{itemize}}
\end{notation}

\newpage

\section{Linear $q$-difference equations of order $1$}\label{sec21}
\subsection{Six fundamental equations}\label{sec211}
In this section, we are interested in the following six equations
\begin{align}
&y(qx)=xy(x)\label{BasicEquation1}\\
&y(qx)=ay(x)\quad\text{with }a\in\C^*\label{BasicEquation2}\\
&y(qx)=\left(1-\dfrac{x}{a}\right)y(x)\quad\text{with }a\in\C^*\label{BasicEquation3}\\
&y(qx)=e^{g(x)}y(x)\quad\text{with }g(x)\in\mathcal O(\C),g(0)=0\label{BasicEquation4}\\
&y(qx)=y(x)+\alpha\quad\text{with }\alpha\in\C^*\label{BasicEquation5}\\
&y(qx)=y(x)+r(x)\quad\text{with }r(x)\in\MerC\text{ analytic at }0\text{ and }r(0)=0\label{BasicEquation6}
\end{align}
For each of them, we make explicit one meromorphic solution on $\C^*$ and we describe the set of all its poles. For the first four equations, we also describe the set of all its zeros.

Before stating our various results (see Propositions \ref{MeromorphicSolutionBasicEquations1}, \ref{MeromorphicSolutionBasicEquations2} and \ref{MeromorphicSolutionBasicEquations3}), let us start by proving the following technical lemma.

\begin{lemma}\label{PreliminariesConvergentPowerSeries}
Let $\displaystyle\sum_{n\geq1}a_nx^n$ be a convergent power series with radius of convergence $0<R\leq+\infty$. Then, the power series
$$\displaystyle\sum_{n\geq1}\dfrac{a_n}{q^n-1}x^n$$
is convergent and its radius of convergence $R'$ is given by $R'=\m{q}R$.
\end{lemma}

\begin{proof}
Lemma \ref{PreliminariesConvergentPowerSeries} is a direct consequence of the Cauchy-Hadamard Theorem and operations on the limits superior. Indeed, we have
$$\underset{n\rightarrow+\infty}{\overline{\lim}}\vert a_n\vert^{1/n}=\dfrac{1}{R}\quad\text{and}\quad\displaystyle\lim_{n\rightarrow+\infty}\vert q^n-1\vert^{1/n}=\vert q\vert;$$
hence, $\underset{n\rightarrow+\infty}{\overline{\lim}}\left\vert\dfrac{a_n}{q^n-1}\right\vert^{1/n}=\dfrac{1}{\vert q\vert R}$.
\end{proof}

Let us now denote by $\Theta_q$ the Jacobi $q$-theta function (see \cite{Ra92}):
$$\Theta_q(x)=\sum_{n\in\Z}(-1)^nq^{-\frac{n(n-1)}{2}}x^n.$$
This function is holomorphic on $\C^*$ and satisfies the Jacobi Triple Product Formula \cite{Jac29}:
$$\Theta_q(x)=\prod_{k\geq0}(1-p^{k+1})(1-xp^k)(1-x^{-1}p^{k+1})\quad\text{with }p=q^{-1}.$$
In particular, its zeros are simple and located at the elements of $q^{\Z}$. Moreover, it satisfies the functional relation
$$\Theta_q(qx)=-qx\Theta_q(x).$$

From these various classical properties, one can easily derive the following result.

\begin{proposition}[\cite{Sau00}]\label{MeromorphicSolutionBasicEquations1}
\begin{enumerate}
\item Equation (\ref{BasicEquation1}) admits the function
$$\frac{\Theta_q^2(x)}{qx\Theta_q(-x)}$$
as solution. It is meromorphic on $\C^*$. Its zeros are double and are located at the elements of $q^{\Z}$; its poles are simple and are located at the elements of $-q^{\Z}$.
\item Equation (\ref{BasicEquation2}) admits the function
$$\frac{\Theta_q(x)}{\Theta_q(a^{-1}x)}$$
as solution. It is constant equal to $1$ if $a=1$ and meromorphic on $\C^*$ otherwise. In the latter case, its zeros and poles are simple and are respectively located at the elements of $q^{\Z}$ and $aq^{\Z}$.
\end{enumerate}
\end{proposition}

In the case of Eqs. (\ref{BasicEquation3}) and (\ref{BasicEquation4}), the situation is much more simpler since the origin $x=0$ is no longer a singular point. In particular, we can make explicit entire solutions as shown in the following.

\begin{proposition}\label{MeromorphicSolutionBasicEquations2}
\begin{enumerate}
\item Equation (\ref{BasicEquation3}) admits the function
$$f_a(x)=\sum_{n\geq0}\frac{1}{(q;q)_n}\left(\frac{x}{a}\right)^n$$
with
$$(q;q)_0=1\quad\text{and}\quad(q;q)_n=\prod_{k=1}^n(1-q^k)$$
as solution. It is entire on $\C$. Its zeros are simple and located at the elements of $aq^{\Na^*}$. Moreover, $f_a(x)=f_1\left(\dfrac{x}{a}\right)$.
\item Equation (\ref{BasicEquation4}) with $g(x)=\displaystyle\sum_{n\geq1}g_nx^n$ admits the function
$$e^{G_g(x)}\quad\text{with}\quad G_g(x)=\sum_{n\geq1}\frac{g_n}{q^n-1}x^n$$
as solution. It is entire on $\C$ and has no zero.
\end{enumerate}
\end{proposition}

\begin{proof}
(1) Looking for the solution $f_a$ in the form $\displaystyle\sum_{n\geq0}a_nx^n$ with $a_0=1$, we get the recurrence relation
$$(q^n-1)a_n = -\dfrac{1}{a}a_{n-1};$$
hence, the identity $a_n=1/((q;q)_na^n)$ for all $n\geq1$.

The function $f_a$ defines obviously an entire function on $\C$ and, from the functional equation
\begin{equation}\label{ProofSolBasicEquation3}
f_a(qx)=\left(1-\dfrac{x}{a}\right)f_a(x),
\end{equation}
it is clear that all the elements of $aq^{\Na^*}$ are zeros of $f_a$. To prove that $f_a$ has no other zero, it is sufficient to observe that if $b\not\in aq^{\Na^*}$ is a zero of $f_a$, then the identity
$$f_a(x)=\left(1-\dfrac{q^{-1}x}{a}\right)...\left(1-\dfrac{q^{-n}x}{a}\right)f_a(q^{-n}x),$$
implies that all the elements $bq^{-n}$ for $n\geq0$ are also zeros of $f_a$, which is impossible. Indeed, the sequence $(bq^{-n})_{n\geq0}$ being convergent to $0$, this would imply $f_a(0)=0$, that is $1=0$.

We are left to prove that the zeros $aq^n$ for $n\geq1$ are simple. Deriving relation (\ref{ProofSolBasicEquation3}) with respect to $x$, we get
$$qf_a'(qx)=-\dfrac{1}{a}f_a(x)+\left(1-\dfrac{x}{a}\right)f'_a(x);$$
hence, the identities
$$\begin{cases}
qf'_a(aq)=-\dfrac{1}{a}f_a(a)\\
qf'_a(aq^{n+1})=(1-q^n)f'_a(aq^n)&\text{for all }n\geq1
\end{cases}.$$
Since $f_a(a)\neq0$, we conclude by recursion on $n$ that $f'_a(aq^n)\neq0$ for all $n\geq1$, which ends the proof of the first point.\\\\
(2) By calculation, we easily check that $G_g(qx)=g(x)+G_g(x)$ and, consequently, that $e^{G_g(x)}$ is a solution of Eq. (\ref{BasicEquation4}). We conclude by observing that Lemma \ref{PreliminariesConvergentPowerSeries} implies that $G_g$ defines an entire function on $\C$.
\end{proof}

For the last two equations (\ref{BasicEquation5}) and (\ref{BasicEquation6}), the situation is more complicated. Indeed, if we can always display a meromorphic solution on $\C^*$, we cannot have a priori only information on its poles.

\begin{proposition}\label{MeromorphicSolutionBasicEquations3}
\begin{enumerate}
\item Equation (\ref{BasicEquation5}) admits the function
$$\alpha x\dfrac{\Theta_q'(x)}{\Theta_q(x)}$$
as solution. It is entire if $\alpha=0$ and meromorphic on $\C^*$ otherwise. In the latter case, its poles are simple and located at the elements of $q^{\Z}$.
\item Equation (\ref{BasicEquation6}) admits the function
$$\dsum{n\geq1}{}{r(p^{n}x)}\quad\text{with }p=q^{-1}$$
as solution. It is meromorphic on $\C$. Its poles are located at the elements of $a q^{\Na^*}$ with order $\mu_{r,a}$ for any $a\in\mathcal P_r$. In particular, it is analytic at the origin and vanishes at $0$.
\end{enumerate}
\end{proposition}

\begin{proof}
(1) The first point is straightforward from the two functional relations $\Theta_q(qx)=-qx\Theta_q(x)$ and $q\Theta_q'(qx)=-q\Theta_q(x)-qx\Theta_q'(x)$.\\\\
(2) To prove the second point, let us first observe that, for all $n\geq1$, the poles of $x\longmapsto r(p^{n}x)$ are the $aq^{n}$'s for any pole $a\in\mathcal P_r$ of $r$. Consequently, setting $\mathcal P=\bigcup_{a\in\mathcal P_r}aq^{\Na^*}$, it is clear that if $x\notin\mathcal P$, then $p^{n}x\notin\mathcal P$ for all $n\geq1$. Observe also that the sequence $(p^{n}x)_{n\geq1}$ tends to $0$ for any $x\in\C$ since $\m{p}<1$.

The function $r$ being analytic at the origin and satisfying $r(0)=0$, there exists a positive constant $C>0$ such that $\m{r(x)}\leq C\m{x}$ for all $\m{x}<\varepsilon$ with a convenient small enough $\varepsilon>0$. From this,  it follows that the series is normally convergent on all the compact sets of $\C\backslash\mathcal P$ (indeed, if $K$ is such a compact set, then there exist $N\geq1$ and $M_K>0$ such that $\m{r(p^nx)}\leq C\m{p^nx}\leq CM_K\m{p}^n$ for all $x\in K$ and all $n\geq N$). The series defines then a meromorphic function on $\C$ with poles in $\mathcal P$, and we can easily check by a direct calculation that it is indeed a solution of Eq. (\ref{BasicEquation6}). Moreover, the order of each of its poles is clearly the one appearing in the statement of Proposition \ref{MeromorphicSolutionBasicEquations3}, (2).
\end{proof}

\subsection{The homogeneous case}\label{sec212}

In this section, we consider a general homogeneous linear $q$-difference equation of the form $y(qx)=m(x)y(x)$ with $m(x)\in\MerCstar$. The case $m\equiv0$ being trivial (the null function is obsviously a solution), we assume in the sequel $m\not\equiv0$.

Since a meromorphic function on $\C^*$ is the quotient of two holomorphic functions on $\C^*$ and since the Birkhoff Decomposition tells us that any holomorphic function on $\C^*$ can be written as a product $h_0(x)h_{\infty}(1/x)$ with two convenient entire functions $h_0(x),h_{\infty}(x)\in\HolC$ (see Proposition \ref{BirkhoffDecomposition} below for more details), the study of this equation is essentially reduced to the case where $m(x)\in\HolC\backslash\{0\}$. For such an equation, Theorem \ref{HomogeneousHolC} below provides an explicit meromorphic solution and describes the set of all its zeros and poles.

Before stating it, let us start by recalling a classical result on homogeneous linear $q$-difference equations which will be very useful to us.

\begin{lemma}\label{CombinationSolutions}
Let $m_1(x),m_2(x)\in\MerCstar$ be two meromorphic functions on $\C^*$ and $H_1(x),H_2(x)\in\MerCstar$ two meromorphic functions on $\C^*$ satisfying the relation
$$H_1(qx)=m_1(x)H_1(x)\quad\text{and}\quad H_2(qx)=m_2(x)H_2(x).$$
Then,
\begin{enumerate}
\item the function $H(x)=H_1(x)H_2(x)\in\MerCstar$ is a meromorphic solution on $\C^*$ of the $q$-difference equation $y(qx)=M(x)y(x)$ with $M(x)=m_1(x)m_2(x)$.
\item the function $\widetilde{H}(x)=H_1(x)/H_2(x)\in\MerCstar$ is a meromorphic solution on $\C^*$ of the $q$-difference equation $y(qx)=\widetilde{M}(x)y(x)$ with $\widetilde{M}(x)=m_1(x)/m_2(x)$.
\end{enumerate}
\end{lemma}

\begin{proof}
By calculations, we have
$$H(qx)=H_1(qx)H_2(qx)=m_1(x)H_1(x)m_2(x)H_2(x)=M(x)H(x);$$
hence, the first point. The second point is proved in a similar way and is left to the reader.
\end{proof}

We are now able to state the main result of this section.

\begin{theorem}\label{HomogeneousHolC}
The equation
\begin{equation}\label{EqHomogeneousHolC}
y(qx)=h(x)y(x)\quad\text{with }h(x)\in\mathcal O(\C)\backslash\{0\}
\end{equation}
admits an entire solution if $h(0)=1$, and a meromorphic solution on $\C^*$ otherwise. Moreover, denoting by $\mu_{h,0}\in\Na$ the order of $0$ as zero of $h$ and setting $\alpha=x^{-\mu_{h,0}}h(x)_{|x=0}$, the zeros and poles of this solution are as stated in Table \ref{ZerosPolesHolomorphicCaseC} below.
\begin{center}
\begin{table}[!h]
\centering\renewcommand{\arraystretch}{1.3}
\begin{tabular}{|p{1cm}|p{4cm}|p{4cm}|}
\cline{2-3}
\multicolumn{1}{c|}{}&zeros located at the elements of...&poles located at the elements of...
\tabularnewline
\hline
$\alpha=1$&$\bullet$ $q^{\Z}$ with order $2\mu_{h,0}$\newline
$\bullet$ $aq^{\Na^*}$ with order $\mu_{h,a}$ for any $a\in\mathcal Z_h^*$&$\bullet$ $-q^{\Z}$ with order $\mu_{h,0}$
\tabularnewline
\hline
$\alpha\neq1$&$\bullet$ $q^{\Z}$ with order $2\mu_{h,0}+1$\newline
$\bullet$ $aq^{\Na^*}$ with order $\mu_{h,a}$ for any $a\in\mathcal Z_h^*$&$\bullet$ $-q^{\Z}$ with order $\mu_{h,0}$\newline
$\bullet$ $\alpha q^{\Z}$ with order $1$
\tabularnewline
\hline
\end{tabular}
\caption{Localization of zeros and poles in the holomorphic case on $\C$}\label{ZerosPolesHolomorphicCaseC}
\end{table}
\end{center}
An explicit writing of such a solution is also given in the constructive proof below.
\end{theorem}

Observe that $\alpha=h(0)$ if and only if $\mu_{h,0}=0$. Observe also that, in the case $\alpha=-1$, the poles of the solution are of order $\mu_{h,0}+1$.

\begin{proof}
The proof is based on the Weierstrass Factorization Theorem. Recall that, since $h\not\equiv0$, the Isolated Zeros Principle implies that the set $\mathcal Z_h$ of all the zeros of $h$ is either empty, or finite, or countable. In all that follows, the constants $\alpha$ and $\mu_{h,0}$ that appear in the decompositions of the function $h$ are those defined in the statement of Theorem \ref{HomogeneousHolC}.\\\\
$\triangleleft$ \textit{First case: $\mathcal Z_h=\emptyset$.} In this case, the function $h$ is written as $h(x)=\alpha e^{g(x)}$ with $g(x)\in\HolC$ an entire function satisfying $g(0)=0$. Applying then Propositions \ref{MeromorphicSolutionBasicEquations1} and \ref{MeromorphicSolutionBasicEquations2} and Lemma \ref{CombinationSolutions}, we deduce that the function
$$\frac{\Theta_q(x)}{\Theta_q(\alpha^{-1}x)}e^{G_g(x)}$$
is a meromorphic solution of (\ref{EqHomogeneousHolC}). More precisely,
\begin{itemize}
\item if $\alpha=1$, it is reduced to $e^{G_g(x)}$; it is therefore entire on $\C$ and has neither zero nor pole;
\item if $\alpha\neq1$, it is meromorphic on $\C^*$; its zeros and poles are simple and are respectively located at the elements of $q^{\Z}$ and $\alpha q^{\Z}$.
\item[]
\end{itemize}
$\triangleleft$ \textit{Second case: $\mathcal Z_h=\{0\}$.} In this case, the function $h$ is written as $h(x)=\alpha x^{\mu_{h,0}} e^{g(x)}$ with $g(x)\in\HolC$ an entire function satisfying $g(0)=0$. Applying again Propositions \ref{MeromorphicSolutionBasicEquations1} and \ref{MeromorphicSolutionBasicEquations2} and Lemma \ref{CombinationSolutions}, we deduce that the function
$$\frac{\Theta_q^{2{\mu_{h,0}}+1}(x)}{q^{\mu_{h,0}}x^{\mu_{h,0}}\Theta_q^{\mu_{h,0}}(-x)\Theta_q(\alpha^{-1}x)}e^{G_g(x)}$$
is a meromorphic solution on $\C^*$ of (\ref{EqHomogeneousHolC}). Moreover,
\begin{itemize}
\item if $\alpha=1$:
\begin{itemize}
\item its zeros are located at the elements of $q^{\Z}$ with order $2\mu_{h,0}$;
\item its poles are located at the elements of $-q^{\Z}$ with order $\mu_{h,0}$.
\end{itemize}
\item if $\alpha\neq1$:
\begin{itemize}
\item its zeros are located at the elements of $q^{\Z}$ with order $2\mu_{h,0}+1$;
\item its poles are located at the elements of $\alpha q^{\Z}$ with order $1$ and at the elements of $-q^{\Z}$ with order $\mu_{h,0}$.
\end{itemize}
\item[]
\end{itemize}
$\triangleleft$ \textit{Third case: $\mathcal Z_h$ finite and $\mathcal Z_h^*\neq\emptyset$.} Let us denote by $a_1$, ..., $a_n$ with $n\geq1$ the nonzero zeros of $h$, and by $v_1$, ..., $v_n$ their respective order of multiplicity. Then, the function $h$ is written as
$$h(x)=\alpha x^{\mu_{h,0}} \left(1-\frac{x}{a_1}\right)^{v_1}...\left(1-\frac{x}{a_n}\right)^{v_n}e^{g(x)}$$
with $g(x)\in\HolC$ an entire function satisfying $g(0)=0$. Applying as before Propositions \ref{MeromorphicSolutionBasicEquations1} and \ref{MeromorphicSolutionBasicEquations2} and Lemma \ref{CombinationSolutions}, we deduce that the function
$$\frac{\Theta_q^{2{\mu_{h,0}}+1}(x)}{q^{\mu_{h,0}}x^{\mu_{h,0}}\Theta_q^{\mu_{h,0}}(-x)\Theta_q(\alpha^{-1}x)}f_{a_1}^{v_1}(x)...f_{a_n}^{v_n}(x)e^{G_g(x)}$$
is a meromorphic solution of (\ref{EqHomogeneousHolC}). More precisely,
\begin{itemize}
\item if $h(0)=1$, then $(\alpha,{\mu_{h,0}})=(1,0)$ and, consequently, the solution is reduced to $f_{a_1}^{v_1}(x)...f_{a_n}^{v_n}(x)e^{G_g(x)}$; it is therefore entire on $\C$ and its zeros are located at the elements of $a_kq^{\Na^*}$ with order $v_k$ for all $k\in\{1,...,n\}$;
\item if $h(0)\neq1$, the solution is meromorphic on $\C^*$ and its zeros and poles are as follows:
\begin{itemize}
\item case $\alpha\neq1$ and ${\mu_{h,0}}=0$: 
\begin{itemize}
\item its zeros are located at the elements of $q^{\Z}$ with order $1$ and at the elements of $a_kq^{\Na^*}$ with order $v_k$ for all $k\in\{1,...,n\}$;
\item its poles are simple and located at the elements of $\alpha q^{\Z}$.
\end{itemize}
\item case $\alpha=1$ and ${\mu_{h,0}}\neq0$: 
\begin{itemize}
\item its zeros are located at the elements of $q^{\Z}$ with order $2{\mu_{h,0}}$ and at the elements of $a_kq^{\Na^*}$ with order $v_k$ for all $k\in\{1,...,n\}$;
\item its poles are located at the elements of $-q^{\Z}$ with order ${\mu_{h,0}}$.
\end{itemize}
\item case $\alpha\neq1$ and ${\mu_{h,0}}\neq0$: 
\begin{itemize}
\item its zeros are located at the elements of $q^{\Z}$ with order $2{\mu_{h,0}}+1$ and at the elements of $a_kq^{\Na^*}$ with order $v_k$ for all $k\in\{1,...,n\}$;
\item its poles are located at the elements of $-q^{\Z}$ with order ${\mu_{h,0}}$, and at the elements of $\alpha q^{\Z}$ with order $1$.
\end{itemize}
\end{itemize}
\item[]
\end{itemize}
$\triangleleft$ \textit{Fourth case: $\mathcal  Z_h$ countable.} We denote by $(a_n)_{n\geq1}$ the set of all the nonzero zeros of $h$, each being counted with its order of multiplicity. According to the Isolated Zeros Principle (recall that $h\not\equiv0$), the sequence $(\m{a_n})_{n\geq1}$ tends to infinity. Then, applying the Weierstrass Factorization Theorem, the function $h$ is written as
$$h(x)=\alpha x^{\mu_{h,0}} e^{g(x)}\prod_{n\geq1}E_{p_n}\left(\frac{x}{a_n}\right)$$
with $g(x)\in\HolC$ an entire function satisfying $g(0)=0$, $(p_n)_{n\geq1}$ a sequence of nonnegative integers such that
\begin{equation}\label{Assumptionpn}
\sum_{n\geq1}\left(\frac{r}{\vert a_n\vert}\right)^{p_n+1}<+\infty\quad\text{for all }r>0,
\end{equation}
and with $E_m(x)$ the Weierstrass' elementary factors defined by
$$E_0(x)=1\text{ and }E_m(x)=(1-x)\exp\left(\sum_{k=1}^{m} \frac{x^k}{k}\right)\text{ for all }m\geq1.$$
Let us now define the entire functions $\widetilde{E}_m(x)$ by
$$\widetilde{E}_0(x)=1\text{ and }\widetilde{E}_m(x)=f_1(x)\exp\left(\sum_{k=1}^{m} \frac{x^k}{k(q^k-1)}\right)\text{ for all }m\geq1,$$
so that $\widetilde{E}_m(qx)=E_m(x)\widetilde{E}_m(x)$ for all $m\geq0$ and all $x\in\C$. Applying then Lemma \ref{ProofHomogeneousHolCLemma1} below, there exist two positive constants $C_1,C_2>0$ such that, for any $r>0$, the following estimate
\begin{equation}\label{EstimateEpn}
\left\vert1-\widetilde{E}_{p_n}\left(\frac{x}{a_n}\right)\right\vert\leq C_1\left(C_2\left\vert\frac{x}{a_n}\right\vert\right)^{p_n+1}\leq C_1\left(\frac{C_2 r}{\vert a_n\vert}\right)^{p_n+1}
\end{equation}
holds for all $n\geq1$ and all $\m{x}\leq r$, as soon as $\m{a_n}\geq r$, inequality valid except at most for a finite number of $n$. Therefore, thanks to the assumption (\ref{Assumptionpn}) on the sequence $(p_n)_{n\geq1}$, we deduce from (\ref{EstimateEpn}) that the series
$$\sum_{n\geq1}\left\vert1-\widetilde{E}_{p_n}\left(\frac{x}{a_n}\right)\right\vert$$
is normally convergent on all the compact sets of $\C$. Consequently, the function $f$ defined by the infinite product
$$f(x)=\prod_{n\geq1}\widetilde{E}_{p_n}\left(\frac{x}{a_n}\right)$$
is entire on $\C$ and satisfies the functional relation
$$f(qx)=\left(\prod_{n\geq1}E_{p_n}\left(\frac{x}{a_n}\right)\right)f(x)$$
for all $x\in\C$. From this, Propositions \ref{MeromorphicSolutionBasicEquations1} and \ref{MeromorphicSolutionBasicEquations2} and Lemma \ref{CombinationSolutions}, we finally derive that the function
$$\frac{\Theta_q^{2{\mu_{h,0}}+1}(x)}{q^{\mu_{h,0}}x^{\mu_{h,0}}\Theta_q^{\mu_{h,0}}(-x)\Theta_q(\alpha^{-1}x)}e^{G_g(x)}\prod_{n\geq1}\widetilde{E}_{p_n}\left(\frac{x}{a_n}\right)$$
is a meromorphic solution of (\ref{EqHomogeneousHolC}). Observing then that
$$\widetilde{E}_{p_n}\left(\frac{x}{a_n}\right)=f_1\left(\frac{x}{a_n}\right)\exp\left(\sum_{k=1}^{p_n} \frac{x^k}{k(q^k-1)a_n^k}\right)=f_{a_n}(x)\exp\left(\sum_{k=1}^{p_n} \frac{x^k}{k(q^k-1)a_n^k}\right)$$
for all $x\in\C$, we get more precisely the following:
\begin{itemize}
\item if $h(0)=1$, then $(\alpha,v)=(1,0)$ and, consequently, the solution is reduced to
$$e^{G_g(x)}\prod_{n\geq1}\widetilde{E}_{p_n}\left(\frac{x}{a_n}\right);$$
it is therefore entire on $\C$ and its zeros are located at the elements of $a_nq^{\Na^*}$ for all $n\geq1$;
\item if $h(0)\neq1$, the solution is meromorphic on $\C^*$ and its zeros and poles are as follows:
\begin{itemize}
\item case $\alpha\neq1$ and ${\mu_{h,0}}=0$: 
\begin{itemize}
\item its zeros are located at the elements of $q^{\Z}$ with order $1$ and at the elements of $a_nq^{\Na^*}$ for all $n\geq1$;
\item its poles are simple and located at the elements of $\alpha q^{\Z}$.
\end{itemize}
\item case $\alpha=1$ and ${\mu_{h,0}}\neq0$: 
\begin{itemize}
\item its zeros are located at the elements of $q^{\Z}$ with order $2{\mu_{h,0}}$ and at the elements of $a_nq^{\Na^*}$ for all $n\geq1$;
\item its poles are located at the elements of $-q^{\Z}$ with order ${\mu_{h,0}}$.
\end{itemize}
\item case $\alpha\neq1$ and ${\mu_{h,0}}\neq0$: 
\begin{itemize}
\item its zeros are located at the elements of $q^{\Z}$ with order $2{\mu_{h,0}}+1$ and at the elements of $a_nq^{\Na^*}$ for all $n\geq1$;
\item its poles are located at the elements of $-q^{\Z}$ with order ${\mu_{h,0}}$, and at the elements of $\alpha q^{\Z}$ with order $1$.
\end{itemize}
\end{itemize}
\end{itemize}
This ends the proof of Theorem \ref{HomogeneousHolC}.
\end{proof}

\begin{lemma}\label{ProofHomogeneousHolCLemma1}
There exist two positive constants $C_1,C_2>0$ such that the following estimate holds for all $m\geq0$ and all $\m{x}\leq1$:
\begin{equation}\label{ProofHomogeneousHolCLemma1Ineq}
\vert1-\widetilde{E}_m(x)\vert\leq C_1(C_2\vert x\vert)^{m+1}.
\end{equation}
\end{lemma}
\begin{proof}
Inequality (\ref{ProofHomogeneousHolCLemma1Ineq}) is valid for any $C_1,C_2>0$ when $m=0$. Let $\varphi_m(x)=1-\widetilde{E}_m(x)$. Then, for all $m\geq1$,
\begin{align*}
\varphi_m'(x)&=-\widetilde{E}_m'(x)\\
&=-\left(f_1'(x)+f_1(x)\sum_{k=0}^{m-1}\frac{x^k}{q^{k+1}-1}\right)\exp\left(\sum_{k=1}^{m} \frac{x^k}{k(q^k-1)}\right)\\
&=-\left(\sum_{n\geq0}\frac{(n+1)x^n}{(q;q)_{n+1}}+\sum_{k=0}^{m-1}\left(\sum_{n\geq k}\frac{x^n}{(q^{k+1}-1)(q;q)_{n-k}}\right)\right)\exp\left(\sum_{k=1}^{m} \frac{x^k}{k(q^k-1)}\right).
\end{align*}
According to the technical Lemma \ref{ProofHomogeneousHolCLemma2} below, all the terms in $x^j$ in the first factor are zero when $j\in\{0,...,m-1\}$. Indeed,
\begin{align*}
\sum_{k=0}^{m-1}\sum_{n=k}^{m-1}\frac{x^n}{(q^{k+1}-1)(q;q)_{n-k}}&=\sum_{n=0}^{m-1}\left(\sum_{k=0}^{n}\frac{1}{(q^{k+1}-1)(q;q)_{n-k}}\right)x^n\\
&=-\dsum{n=0}{m-1}{\frac{(n+1)x^n}{(q;q)_{n+1}}}.
\end{align*}
Therefore, we can write $\varphi_m'(x)$ as
$$\varphi_m'(x)=-\left(\sum_{n\geq m}\frac{(n+1)x^n}{(q;q)_{n+1}}+\sum_{k=0}^{m-1}\left(\sum_{n\geq m}\frac{x^n}{(q^{k+1}-1)(q;q)_{n-k}}\right)\right)\exp\left(\sum_{k=1}^{m} \frac{x^k}{k(q^k-1)}\right).$$
Since all the series which occur define entire functions on $\C$, we derive from this the following estimates for all $m\geq1$ and all $\m{x}\leq1$:
\begin{align*}
\vert\varphi_m'(x)\vert&\leq\vert x\vert^m\left(\sum_{n\geq m}\frac{(n+1)}{(q;q)'_{n+1}}+\sum_{k=0}^{m-1}\left(\sum_{n\geq m}\frac{1}{(\vert q\vert^{k+1}-1)(q;q)'_{n-k}}\right)\right)\exp\left(\sum_{k=1}^{m} \frac{1}{k(\vert q\vert^k-1)}\right)\\
&\leq\vert x\vert^m\left(\sum_{n\geq 0}\frac{(n+1)}{(q;q)'_{n+1}}+\sum_{k=0}^{m-1}\left(\frac{1}{\vert q\vert-1}\sum_{n\geq m-k}\frac{1}{(q;q)'_{n}}\right)\right)\exp\left(\frac{m}{\vert q\vert-1}\right)\\
&\leq\vert x\vert^m\left(\sum_{n\geq 0}\frac{(n+1)}{(q;q)'_{n+1}}+\frac{m}{\vert q\vert-1}\sum_{n\geq 1}\frac{1}{(q;q)'_{n}}\right)\exp\left(\frac{m}{\vert q\vert-1}\right)\\
&\leq m\vert x\vert^m\left(\sum_{n\geq 0}\frac{(n+1)}{(q;q)'_{n+1}}+\frac{1}{\vert q\vert-1}\sum_{n\geq 1}\frac{1}{(q;q)'_{n}}\right)\exp\left(\frac{m}{\vert q\vert-1}\right)\\
&\leq \vert x\vert^m\left(\sum_{n\geq 0}\frac{(n+1)}{(q;q)'_{n+1}}+\frac{1}{\vert q\vert-1}\sum_{n\geq 1}\frac{1}{(q;q)'_{n}}\right)\exp\left(\frac{m\vert q\vert}{\vert q\vert-1}\right),
\end{align*}
where the $(q;q)'_n$ are the positive constants defined by
$$(q;q)'_n=\prod_{k=1}^n(\vert q\vert^k-1)\quad\text{for all }n\geq1.$$
Consequently,
\begin{align*}
\m{1-\widetilde{E}_m(x)}&=\m{\int_{0}^x\varphi_m'(t)dt}\leq C_1(C_2\vert x\vert)^{m+1},
\end{align*}
for all $\m{x}\leq1$, where the positive constants $C_1$ and $C_2$ are respectively defined by
$$C_1=\left(\sum_{n\geq 0}\frac{(n+1)}{(q;q)'_{n+1}}+\frac{1}{\vert q\vert-1}\sum_{n\geq 1}\frac{1}{(q;q)'_{n}}\right)\exp\left(-\frac{\vert q\vert}{\vert q\vert-1}\right)$$
and
$$C_2=\exp\left(\frac{\vert q\vert}{\vert q\vert-1}\right).$$
This completes the proof.
\end{proof}

\begin{lemma}\label{ProofHomogeneousHolCLemma2}
The following identity holds for all $n\geq0$:
$$\frac{n+1}{(q;q)_{n+1}}=-\sum_{k=0}^{n}\frac{1}{(q^{k+1}-1)(q;q)_{n-k}}.$$
\end{lemma}

\begin{proof}
Let us consider the function $\varphi$ defined by
$$\varphi(x)=f_1(x)\exp\left(\sum_{n\geq1}\frac{x^n}{n(q^n-1)}\right)$$
Since $\m{q}>1$ and
$$\displaystyle\sum_{n\geq1}\dfrac{x^n}{n}=-\ln(1-x)$$
for all $\m{x}<1$, we conclude from Lemma \ref{PreliminariesConvergentPowerSeries} that $\varphi(x)$ is well-defined and holomorphic on the disc $\m{x}<1$. Moreover, according to Proposition \ref{MeromorphicSolutionBasicEquations2}, it satisfies the functional relation $\varphi(qx)=\varphi(x)$. Consequently, $\varphi(x)=1$ for all $\m{x}<1$ (indeed, this function is the unique $q$-invariant function which is analytic at the origin). From this, we get $\varphi'(x)=0$ for all $\m{x}<1$; hence, the identity
$$f_1'(x)=-f_1(x)\sum_{n\geq0}\frac{x^n}{q^{n+1}-1}$$
that is
$$\sum_{n\geq0}\frac{(n+1)x^n}{(q;q)_{n+1}}=-\sum_{n\geq0}\left(\sum_{k=0}^n\frac{1}{(q^{k+1}-1)(q;q)_{n-k}}\right)x^n$$
for all $\m{x}<1$, which ends the proof of Lemma \ref{ProofHomogeneousHolCLemma2}.
\end{proof}

As a consequence of Theorem \ref{HomogeneousHolC}, the two following Corollaries \ref{HomogeneousMerC} and \ref{HomogeneousMerCstar} provide explicit meromorphic solutions, as well as the complete description of all their zeros and poles, of the general equation $y(qx)=m(x)y(x)$ with $m(x)\in\MerC\backslash\{0\}$ or $m(x)\in\MerCstar\backslash\{0\}$.

\begin{corollary}\label{HomogeneousMerC}
The equation
\begin{equation}\label{EqHomogeneousMerC}
y(qx)=m(x)y(x)\quad\text{with }m(x)\in\MerC\backslash\{0\}
\end{equation}
admits a meromorphic solution on $\C$ if $0$ is not a pole of $m$ and $m(0)=1$, and a meromorphic solution on $\C^*$ otherwise. Moreover, denoting by
\begin{itemize}
\item$\mu_{m,0}\in\Z$ the order of $0$ at zero (if $\mu_{m,0}\geq0$) or pole (if $\mu_{m,0}<0$) of $m$;
\item$\alpha=x^{-\mu_{m,0}}m(x)_{|x=0}$,
\end{itemize}
the zeros and poles of this solution are as stated in Table \ref{ZerosPolesMeromorphicCaseC1} (case $\mu_{m,0}\geq0$) and in Table \ref{ZerosPolesMeromorphicCaseC2} (case $\mu_{m,0}<0$) below.

\begin{center}
\begin{table}[!h]
\renewcommand{\arraystretch}{1.3}
\begin{tabular}{|p{1cm}|p{4cm}|p{4cm}|}
\cline{2-3}
\multicolumn{1}{c|}{}&zeros located at the elements of...&poles located at the elements of...
\tabularnewline
\hline
$\alpha=1$&$\bullet$ $q^{\Z}$ with order $2\mu_{m,0}$\newline
$\bullet$ $aq^{\Na^*}$ with order $\mu_{m,a}$ for any $a\in\mathcal Z_m^*$&$\bullet$ $-q^{\Z}$ with order $\mu_{m,0}$\newline
$\bullet$ $aq^{\Na^*}$ with order $\mu_{m,a}$ for any $a\in\mathcal P_m^*$
\tabularnewline
\hline
$\alpha\neq1$&$\bullet$ $q^{\Z}$ with order $2\mu_{m,0}+1$\newline
$\bullet$ $aq^{\Na^*}$ with order $\mu_{m,a}$ for any $a\in\mathcal Z_m^*$&$\bullet$ $-q^{\Z}$ with order $\mu_{m,0}$\newline
$\bullet$ $\alpha q^{\Z}$ with order $1$\newline
$\bullet$ $aq^{\Na^*}$ with order $\mu_{m,a}$ for any $a\in\mathcal P_m^*$
\tabularnewline
\hline
\end{tabular}
\caption{Localization of zeros and poles in the meromorphic case on $\C$ with $\mu_{m,0}\geq0$}\label{ZerosPolesMeromorphicCaseC1}
\end{table}
\end{center}
\begin{center}
\begin{table}[!h]
\renewcommand{\arraystretch}{1.3}
\begin{tabular}{|p{1cm}|p{4cm}|p{4cm}|}
\cline{2-3}
\multicolumn{1}{c|}{}&zeros located at the elements of...&poles located at the elements of...
\tabularnewline
\hline
$\alpha=1$&$\bullet$ $-q^{\Z}$ with order $-\mu_{m,0}$\newline
$\bullet$ $aq^{\Na^*}$ with order $\mu_{m,a}$ for any $a\in\mathcal Z_m^*$&$\bullet$ $q^{\Z}$ with order $-2\mu_{m,0}$\newline
$\bullet$ $aq^{\Na^*}$ with order $\mu_{m,a}$ for any $a\in\mathcal P_m^*$
\tabularnewline
\hline
$\alpha\neq1$&$\bullet$ $-q^{\Z}$ with order $-\mu_{m,0}$\newline
$\bullet$ $aq^{\Na^*}$ with order $\mu_{m,a}$ for any $a\in\mathcal Z_m^*$&$\bullet$ $q^{\Z}$ with order $-2\mu_{m,0}-1$\newline
$\bullet$ $\alpha q^{\Z}$ with order $1$\newline
$\bullet$ $aq^{\Na^*}$ with order $\mu_{m,a}$ for any $a\in\mathcal P_m^*$
\tabularnewline
\hline
\end{tabular}
\caption{Localization of zeros and poles in the meromorphic case on $\C$ with $\mu_{m,0}<0$}\label{ZerosPolesMeromorphicCaseC2}
\end{table}
\end{center}
An explicit writing of such a solution can also be obtained by means of Theorem \ref{HomogeneousHolC}.
\end{corollary}

\begin{proof}
Since $m(x)\in\MerC\backslash\{0\}$, then $m$ is written as $m(x)=\alpha x^{\mu_{m,0}} h_1(x)/h_2(x)$ with $h_1(x),h_2(x)\in\HolC$ two entire functions satisfying $h_1(0)=h_2(0)=1$. Doing so, the nonzero zeros (resp. poles) of $m$ are the zeros of $h_1$ (resp. $h_2$). From Theorem \ref{HomogeneousHolC}, there exist two entire functions $H_1(x),H_2(x)\in\HolC$ such that $H_1(qx)=h_1(x)H_1(x)$ and $H_2(qx)=h_2(x)H_2(x)$, the zeros of $H_1$ (resp. $H_2$) being obtained from those of $h_1$ (resp. $h_2$). Applying then Proposition \ref{MeromorphicSolutionBasicEquations1} and Lemma \ref{CombinationSolutions}, we deduce that the function
$$\frac{\Theta_q^{2{\mu_{m,0}}+1}(x)}{q^{\mu_{m,0}}x^{\mu_{m,0}}\Theta_q^{\mu_{m,0}}(-x)\Theta_q(\alpha^{-1}x)}\times\dfrac{H_1(x)}{H_2(x)}$$
is a meromorphic solution of (\ref{EqHomogeneousMerC}) on $\C^*$. The description of the zeros and of the poles of this solution follows from Theorem \ref{HomogeneousHolC}. This ends the proof.
\end{proof}

When $m$ is a general meromorphic function on $\C^*$, the situation is much more complicated and is based on the following multiplicative Birkhoff Decomposition.

\begin{proposition}[Birkhoff Decomposition]\label{BirkhoffDecomposition}
\ \begin{enumerate}
\item Let $h(x)\in\HolCstar\backslash\{0\}$ be a nonzero holomorphic function on $\C^*$.\\
Then, for all $\rho>0$, there exist two constants $\alpha_{\rho}\in\C^*$ and $v_{\rho}\in\Z$, and two entire functions $h_{0,\rho}(x),h_{\infty,\rho}(x)\in\HolC$ satisfying $h_{0,\rho}(0)=h_{\infty,\rho}(0)=1$ such that the following two conditions hold:
\begin{enumerate}
\item the zeros of $h_{0,\rho}(x))$ (resp. $h_{\infty,\rho}(1/x)$) are the elements of $\mathcal Z^*_{h,>\rho}$ (resp. $\mathcal Z^*_{h,\leq\rho}$) with same order of multiplicity;
\item $h(x)=\alpha_{\rho}x^{v_{\rho}}h_{0,\rho}(x)h_{\infty,\rho}(1/x)$ for all $x\in\C^*$.
\end{enumerate}
\item Let $m(x)\in\MerCstar\backslash\{0\}$ be a nonzero meromorphic function on $\C^*$.\\
Then, for all $\rho,\rho'>0$, there exist two constants $\alpha_{\rho,\rho'}\in\C^*$ and $v_{\rho,\rho'}\in\Z$, and two meromorphic functions $m_{0,\rho,\rho'}(x),m_{\infty,\rho,\rho'}(x)\in\MerC$ without pole at $0$ and satisfying $m_{0,\rho,\rho'}(0)=m_{\infty,\rho,\rho'}(0)=1$ such that the following three conditions hold:
\begin{enumerate}
\item the zeros of $m_{0,\rho,\rho'}(x))$ (resp. $m_{\infty,\rho,\rho'}(1/x)$) are the elements of $\mathcal Z^*_{m,>\rho}$ (resp. $\mathcal Z^*_{m,\leq\rho}$) with same order of multiplicity;
\item the poles of $m_{0,\rho,\rho'}(x))$ (resp. $m_{\infty,\rho,\rho'}(1/x)$) are the elements of $\mathcal P^*_{m,>\rho'}$ (resp. $\mathcal P^*_{m,\leq\rho'}$) with same order of multiplicity;
\item $m(x)=\alpha_{\rho,\rho'}x^{v_{\rho,\rho'}}m_{0,\rho,\rho'}(x)m_{\infty,\rho,\rho'}(1/x)$ for all $x\in\C^*$.
\end{enumerate}
\end{enumerate}
\end{proposition}

Observe that the distribution of the initial zeros and poles of $h(x)$ and $m(x)$ in the two decompositions above are totally arbitrary and are therefore left to a free choice.

\begin{proof}
Since $m(x)\in\MerCstar\backslash\{0\}$ is the quotient of two nonzero holomorphic functions on $\C^*$, it is sufficient to prove the first point. We reproduce below the proof provided to us by A. Ancona and whom we wish to thank warmly here.

So, let us consider a holomorphic function $h(x)\in\HolCstar\backslash\{0\}$ and let us start by considering the set $\mathcal Z_h^*$ of all its nonzero zeros.
\begin{itemize}
\item When $\mathcal Z^*_{h,>\rho}\neq\emptyset$, its elements have no accumulation point in $\C$ and, from the Weiestrass Theorem (see for instance \cite[Th\'{e}or\`{e}me 15.9 page 282]{Rud78}), there exists an entire function $f_{0,\rho}(x)\in\HolC$ whose zeros are the elements of $\mathcal Z^*_{h,>\rho}$ with same order of multiplicity.\\
When $\mathcal Z^*_{h,>\rho}=\emptyset$, we choose for $f_{0,\rho}$ the function $f_{0,\rho}\equiv1$.
\item When $\mathcal Z^*_{h,\leq\rho}\neq\emptyset$, its elements $a$ may have $0$ as accumulation point (and it is the only one possible!). Therefore, the set of all their inverse $1/a$ has no accumulation point in $\C$ and, applying again the Weierstrass Theorem, there exists an entire function $f_{\infty,\rho}(x)\in\HolC$ whose zeros are the elements $1/a$ for any $a\in\mathcal Z^*_{h,\leq\rho}$ with same order of multiplicity than $a$. In particular, $f_{\infty,\rho}(1/x)\in\HolCstar$ and its zeros are the elements of $\mathcal Z^*_{h,\leq\rho}$ with same order of multiplicity.\\
When $\mathcal Z^*_{h,\leq\rho}=\emptyset$, we choose for $f_{\infty,\rho}$ the function $f_{\infty,\rho}\equiv1$.
\end{itemize}
According to our assumptions on $f_{0,\rho}$ and $f_{\infty,\rho}$, the function $k_{\rho}$ defined by
$$k_{\rho}(x)=\dfrac{h(x)}{f_{0,\rho}(x)f_{\infty,\rho}(1/x)}$$
is holomorphic on $\C^*$ and without zero in $\C^*$. Therefore, the function $k_{\rho}(e^x)$ is entire without zero in $\C$ and, consequently, there exists an entire function $\varphi_{\rho}(x)\in\HolC$ such that
$$k_{\rho}(e^x)=e^{\varphi_{\rho}(x)}\quad\text{for all }x\in\C.$$
Let us now observe that $x\longmapsto k_{\rho}(e^x)$ is $2i\pi$-periodic; hence, $\varphi_{\rho}(x+2i\pi)-\varphi_{\rho}(x)\in2i\pi\Z$ for all $x\in\C$. Since the function $x\longmapsto\varphi_{\rho}(x+2i\pi)-\varphi_{\rho}(x)$ is also continuous, we deduce there exists an integer $v_{\rho}\in\Z$ such that
$$\varphi_{\rho}(x+2i\pi)-\varphi_{\rho}(x)=2i\pi v_{\rho}\quad\text{for all }x\in\C.$$
Therefore, the function $\psi_{\rho}(x)=\varphi_{\rho}(x)-v_{\rho}x$ is $2i\pi$-periodic and we have
$$k_{\rho}(e^x)=e^{v_{\rho}x}e^{\psi_{\rho}(x)}\quad\text{for all }x\in\C.$$
Going back to the function $k$ itself, we get the identity
$$k_{\rho}(x)=x^{v_{\rho}}e^{\psi_{\rho}(\ln(x))}\quad\text{for all }x\in\C^*,$$
where, accordingly the $2i\pi$-periodicity of $\psi_{\rho}$, the composition $x\longmapsto\psi_{\rho}(\ln(x))$ is holomorphic univalent on $\C^*$. In particular, it can be decomposed into a Laurent series at $0$: there exist two entire functions $\psi_{0,\rho}(x),\psi_{\infty,\rho}(x)\in\HolC$ such that
$$\psi_{\rho}(\ln(x))=\psi_{0,\rho}(x)+\psi_{\infty,\rho}\left(\dfrac{1}{x}\right)\quad\text{for all }x\in\C^*.$$
Hence, the identity
$$h(x)=x^{v_{\rho}}f_{0,\rho}(x)f_{\infty,\rho}\left(\dfrac{1}{x}\right)e^{\psi_{0,\rho}(x)}e^{\psi_{\infty,\rho}(1/x)}$$
for all $x\in\C^*$. The choices
\begin{align*}
&h_{0,\rho}(x)=\dfrac{f_{0,\rho}(x)e^{\psi_{0,\rho}(x)}}{f_{0,\rho}(0)e^{\psi_{0,\rho}(0)}}\in\HolC\\
&h_{\infty,\rho}(x)=\dfrac{f_{\infty,\rho}(x)e^{\psi_{\infty,\rho}(x)}}{f_{\infty,\rho}(0)e^{\psi_{\infty,\rho}(0)}}\in\HolC\\
&\alpha_{\rho}=f_{0,\rho}(0)f_{\infty,\rho}(0)e^{\psi_{0,\rho}(0)}e^{\psi_{\infty,\rho}(0)}
\end{align*}
complete the proof.
\end{proof}

\begin{corollary}\label{HomogeneousMerCstar}
Let $\rho,\rho'>0$ be two positive real numbers and $m(x)\in\MerCstar\backslash\{0\}$ a meromorphic function on $\C^*$ written in the form
$$m(x)=\alpha_{\rho,\rho'}x^{v_{\rho,\rho'}}m_{0,\rho,\rho'}(x)m_{\infty,\rho,\rho'}(1/x)$$
as in Proposition \ref{BirkhoffDecomposition}. Then, the equation
\begin{equation}\label{EqHomogeneousMerCstar}
y(qx)=m(x)y(x)
\end{equation}
admits a meromorphic solution on $\C^*$ whose the zeros and poles are as stated in Table \ref{ZerosPolesMeromorphicCaseCstar1} (case $v_{\rho,\rho'}\geq0$) and in Table \ref{ZerosPolesMeromorphicCaseCstar2} (case $v_{\rho,\rho'}<0$) below.
\begin{center}
\begin{table}[!h]
\renewcommand{\arraystretch}{1.3}
\begin{tabular}{|p{1.5cm}|p{4cm}|p{4cm}|}
\cline{2-3}
\multicolumn{1}{c|}{}&zeros located at the elements of...&poles located at the elements of...
\tabularnewline
\hline
$\alpha_{\rho,\rho'}=1$&$\bullet$ $q^{\Z}$ with order $2v_{\rho,\rho'}$\newline
$\bullet$ $aq^{\Na^*}$ with order $\mu_{m,a}$ for any $a\in\mathcal Z_{m,>\rho}^*$\newline
$\bullet$ $aq^{-\Na}$ with order $\mu_{m,a}$ for any $a\in\mathcal P_{m,\leq\rho'}^*$&$\bullet$ $-q^{\Z}$ with order $v_{\rho,\rho'}$\newline
$\bullet$ $aq^{\Na^*}$ with order $\mu_{m,a}$ for any $a\in\mathcal P_{m,>\rho'}^*$\newline
$\bullet$ $aq^{-\Na}$ with order $\mu_{m,a}$ for any $a\in\mathcal Z_{m,\leq\rho}^*$
\tabularnewline
\hline
$\alpha_{\rho,\rho'}\neq1$&$\bullet$ $q^{\Z}$ with order $2v_{\rho,\rho'}+1$\newline
$\bullet$ $aq^{\Na^*}$ with order $\mu_{m,a}$ for any $a\in\mathcal Z_{m,>\rho}^*$\newline
$\bullet$ $aq^{-\Na}$ with order $\mu_{m,a}$ for any $a\in\mathcal P_{m,\leq\rho'}^*$&$\bullet$ $-q^{\Z}$ with order $v_{\rho,\rho'}$\newline
$\bullet$ $\alpha_{\rho,\rho'} q^{\Z}$ with order $1$\newline
$\bullet$ $aq^{\Na^*}$ with order $\mu_{m,a}$ for any $a\in\mathcal P_{m,>\rho'}^*$\newline
$\bullet$ $aq^{-\Na}$ with order $\mu_{m,a}$ for any $a\in\mathcal Z_{m,\leq\rho}^*$
\tabularnewline
\hline
\end{tabular}
\caption{Localization of zeros and poles in the meromorphic case on $\C^*$ with $v_{\rho,\rho'}\geq0$}\label{ZerosPolesMeromorphicCaseCstar1}
\end{table}
\end{center}
\begin{center}
\begin{table}[!h]
\renewcommand{\arraystretch}{1.3}
\begin{tabular}{|p{1.5cm}|p{4cm}|p{4cm}|}
\cline{2-3}
\multicolumn{1}{c|}{}&zeros located at the elements of...&poles located at the elements of...
\tabularnewline
\hline
$\alpha_{\rho,\rho'}=1$&$\bullet$ $-q^{\Z}$ with order $-v_{\rho,\rho'}$\newline
$\bullet$ $aq^{\Na^*}$ with order $\mu_{m,a}$ for any $a\in\mathcal Z_{m,>\rho}^*$\newline
$\bullet$ $aq^{-\Na}$ with order $\mu_{m,a}$ for any $a\in\mathcal P_{m,\leq\rho'}^*$&$\bullet$ $q^{\Z}$ with order $-2v_{\rho,\rho'}$\newline
$\bullet$ $aq^{\Na^*}$ with order $\mu_{m,a}$ for any $a\in\mathcal P_{m,>\rho'}^*$\newline
$\bullet$ $aq^{-\Na}$ with order $\mu_{m,a}$ for any $a\in\mathcal Z_{m,\leq\rho}^*$
\tabularnewline
\hline
$\alpha_{\rho,\rho'}\neq1$&$\bullet$ $-q^{\Z}$ with order $-v_{\rho,\rho'}$\newline
$\bullet$ $aq^{\Na^*}$ with order $\mu_{m,a}$ for any $a\in\mathcal Z_{m,>\rho}^*$\newline
$\bullet$ $aq^{-\Na}$ with order $\mu_{m,a}$ for any $a\in\mathcal P_{m,\leq\rho'}^*$&$\bullet$ $q^{\Z}$ with order $-2v_{\rho,\rho'}-1$\newline
$\bullet$ $\alpha_{\rho,\rho'} q^{\Z}$ with order $1$\newline
$\bullet$ $aq^{\Na^*}$ with order $\mu_{m,a}$ for any $a\in\mathcal P_{m,>\rho'}^*$\newline
$\bullet$ $aq^{-\Na}$ with order $\mu_{m,a}$ for any $a\in\mathcal Z_{m,\leq\rho}^*$
\tabularnewline
\hline
\end{tabular}
\caption{Localization of zeros and poles in the meromorphic case on $\C^*$ with $v_{\rho,\rho'}<0$}\label{ZerosPolesMeromorphicCaseCstar2}
\end{table}
\end{center}
\end{corollary}

\begin{proof}
Applying Corollary \ref{HomogeneousMerC}, there exist two meromorphic functions
$$M_{0,\rho,\rho'}(x),M_{\infty,\rho,\rho'}(x)\in\MerC$$ such that
\begin{align*}
&M_{0,\rho,\rho'}(qx)=m_{0,\rho,\rho'}(x)M_{0,\rho,\rho'}(x)\text{ and}\\
&M_{\infty,\rho,\rho'}(qx)=m_{\infty,\rho,\rho'}^{-1}(x)M_{\infty,\rho,\rho'}(x).
\end{align*}
Then, according to Proposition \ref{MeromorphicSolutionBasicEquations1} and Lemma \ref{CombinationSolutions}, the function
$$\frac{\Theta_q^{2v_{\rho,\rho'}+1}(x)}{q^{v_{\rho,\rho'}}x^{v_{\rho,\rho'}}\Theta_q^{v_{\rho,\rho'}}(-x)\Theta_q(\alpha_{\rho,\rho'}^{-1}x)}\times M_{0,\rho,\rho'}(x)M_{\infty,\rho,\rho'}\left(\frac{q}{x}\right)$$
is a meromorphic solution of (\ref{EqHomogeneousMerCstar}) on $\C^*$, and the description of the zeros and poles follows from Corollary \ref{HomogeneousMerC}.
\end{proof}

We end this section with some examples.

\begin{example}
\textnormal{\begin{enumerate}
\item As a first example, let us consider the equation
\begin{equation}\label{HomEqEx0}
y(qx)=(1-x^3)y(x),
\end{equation}
with a polynomial coefficient. Since $(1-x^3)_{|x=0}=1$ and since its zeros are simple and located at the cubic roots $1$, $j=-1/2+i\sqrt{3}/2$ and $j^2$ of the unit, we easily derived from Theorem \ref{HomogeneousHolC} that Eq. (\ref{HomEqEx0}) admits an entire solution. More precisely, observing that $1-x^3=(1-x)(1-x/j)(1-x/j^2)$, it is given by the function
$$f_1(x)f_j(x)f_{j^2}(x)$$
with $f_a(x)$ as in Proposition \ref{MeromorphicSolutionBasicEquations2}, and its zeros are simple and located at the elements of $q^{\Na^*}$, $jq^{\Na^*}$ and $j^2q^{\Na^*}$.
\item Let us then consider the equation
\begin{equation}\label{HomEqEx1}
y(qx)=\sin(x)y(x).
\end{equation}
From the Weierstrass Factorization Theorem, we have
$$\sin(x)=x\prod_{n\geq1}\left(1-\dfrac{x}{n\pi}\right)\left(1+\dfrac{x}{n\pi}\right)\quad\text{for all }x\in\C.$$
Applying then Theorem \ref{HomogeneousHolC}, we deduce that a meromorphic solution on $\C^*$ of Eq. (\ref{HomEqEx1}) is given by the function
$$\dfrac{\Theta_q^2(x)}{qx\Theta_q(-x)}\prod_{n\geq1}f_{n\pi}(x)f_{-n\pi}(x).$$
Moreover,
\begin{itemize}
\item its zeros are located at the elements of $q^{\Z}$ with order $2$ and at the elements of $n\pi q^{\Na^*}$ with order $1$ for all $n\in\Z\backslash\{0\}$;
\item its poles are simple and located at the elements of $-q^{\Z}$.
\end{itemize}
\item Let us now consider the equation
\begin{equation}\label{HomEqEx2}
y(qx)=\Gamma(x)y(x).
\end{equation}
The function $\Gamma$ is meromorphic on $\C$ and, from the Weierstrass Factorization Theorem, we have
$$\Gamma(x)=\dfrac{e^{-\gamma x}}{x\displaystyle\prod_{n\geq1}\left(1+\dfrac{x}{n}\right)e^{-x/n}}\quad\text{for all }x\in\C\backslash(-\Na).$$
Consequently, Theorem \ref{HomogeneousHolC} and Corollary \ref{HomogeneousMerC} tell us that a meromorphic solution on $\C^*$ of Eq. (\ref{HomEqEx2}) is given by the function
$$\dfrac{qxe^{\gamma x/(1-q)}\Theta_q(-x)}{\Theta^2_q(x)\displaystyle\prod_{n\geq1}f_{-n}(x)e^{x/(n(1-q))}}.$$
Moreover,
\begin{itemize}
\item its zeros are simple and located at the elements of $-q^{\Z}$;
\item its poles are located at the elements of $q^{\Z}$ with order $2$ and at the elements of $-n q^{\Na^*}$ with order $1$ for all integer $n\geq1$.
\end{itemize}
\item As a final example, let us consider the equation
\begin{equation}\label{HomEqEx3}
y(qx)=\sin\left(\dfrac{2}{x}\right)y(x).
\end{equation}
We have $\sin(2/x)\in\HolCstar\subset\MerCstar$ and, for its Birkhoff decomposition, we choose the one provided by the Weierstrass Factorization Theorem:
$$\sin\left(\dfrac{2}{x}\right)=2x^{-1}h_\infty\left(\dfrac{1}{x}\right)$$
where $h_\infty(x)$ is the entire function defined by
$$h_\infty(x)=\prod_{n\geq1}\left(1-\dfrac{2x}{n\pi}\right)\left(1+\dfrac{2x}{n\pi}\right)\quad\text{for all }x\in\C.$$
In particular, the zeros of $h_{\infty}(1/x)$ are the nonzero zeros of $\sin(2/x)$: they are simple and located at the elements $2/(n\pi)$ with $n\in\Z^*$. Using then Corollary \ref{HomogeneousMerCstar}, Table \ref{ZerosPolesMeromorphicCaseCstar2}, we deduce that Eq. (\ref{HomEqEx3}) admits a meromorphic solution on $\C^*$, whose zeros and poles are as follows:
\begin{itemize}
\item its zeros are simple and located at the elements of $-q^{\Z}$;
\item its poles are located at the elements of $q^{\Z}$ with order $3$, at the elements of $2q^{\Z}$ with order $1$, and at the elements of $\dfrac{2}{n\pi}q^{-\Na}$ with order $1$ for all integer $n\in\Z^*$.
\end{itemize}
\end{enumerate}}
\end{example}

\subsection{The inhomogeneous case}\label{SectionInhomogeneousLinearQDiffOrder1}
In this section, we are interested in the general inhomogeneous linear $q$-difference equation
\begin{equation}\label{EqInhomogeneousMerCstar}
y(qx)=m(x)y(x)+r(x)\quad\text{with }m(x),r(x)\in\MerCstar.
\end{equation}
The cases $m\equiv0$ and $r\equiv0$ being trivial (the first one provides the meromorphic solution $r(x/q)$ and the second one coincides with the linear case), we assume in the sequel $m\not\equiv0$ and $r\not\equiv0$.

From Corollary \ref{HomogeneousMerCstar}, there exists a meromorphic solution $M(x)\in\MerCstar$ of the equation $y(qx)=m(x)y(x)$. Applying then the change of unknown function $y(x)=M(x)z(x)$, Eq. (\ref{EqInhomogeneousMerCstar}) becomes
\begin{equation}\label{TransformedEqInhomogeneousMerCstar}
z(qx)=z(x)+R(x)\quad\text{with }R(x)=M^{-1}(qx)r(x)\in\MerCstar.
\end{equation}

Proposition \ref{InhomogeneousRMerCstar} below provides a meromorphic solution on $\C^*$ of this equation, as well as the complete description of all its poles. The construction of such a solution is based on the following additive decomposition of $R(x)$ which generalizes the Laurent series decomposition.

\begin{proposition}\label{AdditiveDecompositonMCstar}
Let $R(x)\in\MerCstar$ be a meromorphic function on $\C^*$.\\
Then, for all $\rho >0$, there exist a constant $\alpha_{\rho}\in\C$ and two meromorphic functions $R_{0,\rho}(x),R_{\infty,\rho}(x)\in\MerC$ without pole at $0$ and satisfying $R_{0,\rho}(0)=R_{\infty,\rho}(0)=0$ such that the following three conditions hold:
\begin{enumerate}
\item the poles of $R_{0,\rho}(x)$ (resp. $R_{\infty,\rho}(1/x)$) are the elements of $\mathcal P^*_{R,>\rho}$ (resp. $\mathcal P^*_{R,\leq\rho}$);
\item for each $a\in\mathcal P^*_{R,>\rho}$ (resp. $\mathcal P^*_{R,\leq\rho}$), the principal part of $R_{0,\rho}(x)$ (resp. $R_{\infty,\rho}(1/x)$) at $a$ coincides with the principal part of $R(x)$ at $a$;
\item$R(x)=R_{0,\rho}(x)+\alpha_{\rho}+R_{\infty,\rho}(1/x)$ for all $x\in\C^*$.
\end{enumerate}
\end{proposition}

\begin{proof}
When $\mathcal P^*_R=\emptyset$, the function $R$ is holomorphic on $\C^*$ and the decomposition stems obvious from the decomposition of $R$ into Laurent series at $0$. In particular, the functions $R_{0,\rho}(x)$ and $R_{\infty,\rho}(x)$ are entire.\\
Let us now suppose $\mathcal P^*_R\neq\emptyset$.
\begin{itemize}
\item When $\mathcal P^*_{R,>\rho}\neq\emptyset$, its elements have no accumulation point in $\C$ and, from the Mittag-Leffler Theorem \cite{MitLef1884} (see also \cite[Th\'{e}or\`{e}me 15.13 page 285]{Rud78}), there exists a meromorphic function $f_{0,\rho}(x)\in\MerC$ whose poles are the elements of $\mathcal P^*_{R,>\rho}$ and whose principal part at each $a\in\mathcal P^*_{R,>\rho}$ is $P_{R,a}(x)$. In particular, this function is analytic at the origin.\\
When $\mathcal P^*_{R,>\rho}=\emptyset$, we choose for $f_{0,\rho}$ the null function.
\item When $\mathcal P^*_{R,\leq\rho}\neq\emptyset$, the previous reasoning does not apply anymore because $0$ may be an accumulation point (and it is the only one possible!). To get around this difficulty, we therefore consider, not the set of elements $a\in\mathcal P^*_{R,\leq\rho}$, but the set of their inverse $1/a$ (which is well-defined since all the poles of $R(x)$ are nonzero). By construction, this new set has no accumulation point in $\C$. For any $a\in\mathcal P^*_{R,\leq\rho}$, we denote by $P_{1/a}(x)$ the unique polynomial in $1/(x-1/a)$ without constant term such that the principal part of $P_{1/a}(1/x)$ at $a$, that is the polynomial in $1/(x-a)$ without constant term, is $P_{R,a}(x)$. Then, applying again the Mittag-Leffler Theorem, there exists a meromorphic function $f_{\infty,\rho}(x)\in\MerC$ whose poles are the elements $1/a$ for any $a\in\mathcal P^*_{R,\leq\rho}$ and whose principal part at each point $1/a$ is $P_{1/a}(x)$. In particular, this function is analytic at the origin. Moreover, the function $f_{\infty,\rho}(1/x)$ is meromorphic on $\C^*$, its poles are the elements of $\mathcal P^*_{R,\leq\rho}$ and the principal part at each $a\in\mathcal P^*_{R,\leq\rho}$ is $P_{R,a}(x)$.\\
When $\mathcal P^*_{R,\leq\rho}=\emptyset$, we choose for $f_{\infty,\rho}$ the null function.
\end{itemize}
According to our assumptions on $f_{0,\rho}$ and $f_{\infty,\rho}$, the function $R(x)-f_{0,\rho}(x)-f_{\infty,\rho}(1/x)$ is holomorphic on $\C^*$ and can then be decomposed into a Laurent series at $0$: there exist two entire functions $g_0(x),g_{\infty}(x)\in\HolC$ such that
$$R(x)-f_{0,\rho}(x)-f_{\infty,\rho}(1/x)=g_0(x)+g_{\infty}(1/x)$$
for all $x\in\C^*$. The choices
\begin{align*}
&R_{0,\rho}(x)=f_{0,\rho}(x)+g_0(x)-f_{0,\rho}(0)-g_0(0)\in\MerC\\
&R_{\infty,\rho}(x)=f_{\infty,\rho}(x)+g_{\infty}(x)-f_{\infty,\rho}(0)-g_{\infty}(0)\in\MerC\\
&\alpha_{\rho}=f_{0,\rho}(0)+g_0(0)+f_{\infty,\rho}(0)+g_{\infty}(0)
\end{align*}
complete the proof.
\end{proof}

\begin{remark}
\textnormal{\begin{itemize}
\item Unlike the Laurent series decomposition, the decomposition obtained in Proposition \ref{AdditiveDecompositonMCstar} is not unique, even for a fixed $\rho>0$, since it depends on the choice of the two meromorphic functions $f_{0,\rho}(x)$ and $f_{\infty,\rho}(x)$.
\item As in the case of the Birkhoff Decomposition (see Proposition \ref{BirkhoffDecomposition}), the distribution of the initial poles of $R(x)$ in the above decomposition is still completely arbitrary and is therefore left to a free choice.
\end{itemize}}
\end{remark}

We are now able to solve Eq. (\ref{TransformedEqInhomogeneousMerCstar}).

\begin{proposition}\label{InhomogeneousRMerCstar}
For all $\rho>0$, the equation
\begin{equation}\label{EqInhomogeneousRMerCstar}
y(qx)=y(x)+R(x)\quad\text{with }R(x)\in\MerCstar
\end{equation}
admits a meromorphic solution on $\C^*$ whose the poles are located at the elements of
\begin{itemize}
\item$q^{\Z}$ with order at most $1$;
\item$aq^{\Na^*}$with order $\mu_{R,a}$ for any $a\in\mathcal P^*_{R,>\rho}$;
\item$aq^{-\Na}$with order $\mu_{R,a}$ for any $a\in\mathcal P^*_{R,\leq\rho}$.
\end{itemize}
Moreover, an explicit writing of such a solution is given in the constructive proof below.
\end{proposition}

\begin{proof}
Let us write $R(x)$ in the form $R(x)=R_{0,\rho}(x)+\alpha_{\rho}+R_{\infty,\rho}(1/x)$ with $\alpha_{\rho}\in\C$ and $R_{0,\rho}(x),R_{\infty,\rho}(x)\in\MerC$ as in Proposition \ref{AdditiveDecompositonMCstar}. Applying Proposition \ref{MeromorphicSolutionBasicEquations3}, we get the following solutions:
\begin{itemize}
\item the equation $z(qx)=z(x)+\alpha_{\rho}$ admits the function
$$z_{\alpha_{\rho}}(x)=\alpha_{\rho} x\dfrac{\Theta_q'(x)}{\Theta_q(x)}$$
as solution. It is entire if $\alpha_{\rho}=0$ and meromorphic on $\C^*$ otherwise. In the latter case, its pole are simple and located at the elements of $q^{\Z}$;
\item the equation $z(qx)=z(x)+R_{0,\rho}(x)$ admits the function
$$z_{0,\rho}(x)=\dsum{n\geq1}{}{R_{0,\rho}(p^nx)}\quad\text{with }p=q^{-1}$$
as solution. It is meromorphic on $\C$ and its poles are located at the elements of $aq^{\Na^*}$ with order $\mu_{R,a}$ for any $a\in\mathcal P^*_{R,>\rho}$;
\item the equation $z(qx)=z(x)-R_{\infty,\rho}(x)$ admits the function
$$z_{\infty,\rho}(x)=-\dsum{n\geq1}{}{R_{\infty,\rho}(p^nx)}$$
as solution. It is meromorphic on $\C$ and its poles are located at the elements of $\dfrac{1}{a}q^{\Na^*}$ with order $\mu_{R,a}$ for any $a\in\mathcal P^*_{R,\leq\rho}$.
\end{itemize}
Since $z_{\infty,\rho}(q/x)$ is a meromorphic solution on $\C^*$ of the equation
$$z(qx)=z(x)+R_{\infty,\rho}(1/x)$$
with poles located at the elements of $aq^{-\Na}$ with order $\mu_{R,a}$ for any $a\in\mathcal P^*_{R,\leq\rho}$, we deduce from the Superposition Principle that the function
$$z_{0,\rho}(x)+z_{\alpha_{\rho}}(x)+z_{\infty,\rho}(q/x)=\dsum{n\geq1}{}{R_{0,\rho}(p^nx)}+\alpha_{\rho} x\dfrac{\Theta_q'(x)}{\Theta_q(x)}-\dsum{n\geq0}{}{R_{\infty,\rho}\left(\dfrac{p^n}{x}\right)}$$
is a meromorphic solution on $\C^*$ of Eq. (\ref{EqInhomogeneousRMerCstar}). The full description of its poles follows from the previous ones.
\end{proof}

\begin{corollary}\label{InhomogeneousMerCstar}
Let $\rho,\rho',\rho''>0$. Then, the equation
\begin{equation}\label{EqInhomogeneousMerCstarNonzero}
y(qx)=m(x)y(x)+r(x)\quad\text{with }m(x),r(x)\in\MerCstar\backslash\{0\}
\end{equation}
admits a meromorphic solution on $\C^*$ whose poles are located at:
\begin{enumerate}
\item Case $m(x)\in\MerC\backslash\{0\}$:
\begin{itemize}
\item $\pm q^{\Z}$ and $\alpha q^{\Z}$;
\item$aq^{\Na^*}$ for any $a\in\mathcal P_m^*\cup\mathcal P_{>\rho}$;
\item$aq^{-\Na}$ for any $a\in\mathcal P_{\leq\rho}$,
\end{itemize}
where we set
\begin{itemize}
\item$\alpha=x^{-v}m(x)_{|x=0}$ with $v$ the order of $0$ at zero/pole of $m$;
\item$\mathcal P=\mathcal P_r^*\cup\{aq^{\Na};a\in\mathcal Z_m^*\}$;
\item$\mathcal P_{\leq\rho}=\{a\in\mathcal P;\m{a}\leq\rho\}$;
\item$\mathcal P_{>\rho}=\{a\in\mathcal P;\m{a}>\rho\}$.
\end{itemize}
\item  Case $m(x)\in\MerCstar\backslash\{0\}$:
\begin{itemize}
\item $\pm q^{\Z}$ and $\alpha_{\rho',\rho''} q^{\Z}$;
\item$aq^{\Na^*}$ for any $a\in\mathcal P_{m,>\rho''}^*\cup\mathcal P_{>\rho}$;
\item$aq^{-\Na}$ for any $a\in\mathcal Z_{m,\leq\rho'}^*\cup\mathcal P_{\leq\rho}$,
\end{itemize}
where we set
\begin{itemize}
\item$\alpha_{\rho',\rho''}$ as in Corollary \ref{HomogeneousMerCstar};
\item$\mathcal P=\mathcal P_r^*\cup\{aq^{\Na};a\in\mathcal Z_{m,>\rho'}^*\}\cup\{aq^{-\Na^*};a\in\mathcal P_{m,\leq\rho''}^*\}$;
\item$\mathcal P_{\leq\rho}=\{a\in\mathcal P;\m{a}\leq\rho\}$;
\item$\mathcal P_{>\rho}=\{a\in\mathcal P;\m{a}>\rho\}$.
\end{itemize}
\end{enumerate}
\end{corollary}

\begin{proof}
From the calculations made at the beginning of Section \ref{SectionInhomogeneousLinearQDiffOrder1}, we deduce from Proposition \ref{InhomogeneousRMerCstar} that a meromorphic solution on $\C^*$ of Eq. (\ref{EqInhomogeneousMerCstarNonzero}) is given by the function
\begin{equation}\label{SolEqInhomogeneousMerCstarNonzero}
M(x)\left(\dsum{n\geq1}{}{R_{0,\rho}(p^nx)}+\alpha_{\rho} x\dfrac{\Theta_q'(x)}{\Theta_q(x)}-\dsum{n\geq0}{}{R_{\infty,\rho}\left(\dfrac{p^n}{x}\right)}\right),
\end{equation}
where
\begin{itemize}
\item$M(x)$ is a meromorphic solution on $\C^*$ of the equation $y(qx)=m(x)y(x)$ as in Corollary \ref{HomogeneousMerC} when $m(x)\in\MerC$ and as in Corollary \ref{HomogeneousMerCstar} when $m(x)\in\MerCstar$;
\item$\alpha_{\rho}\in\C$ and $R_{0,\rho}(x),R_{\infty,\rho}(x)\in\MerC$ are the elements of the decomposition of $R(x)=M^{-1}(qx)r(x)$ given by Proposition \ref{AdditiveDecompositonMCstar}.
\end{itemize}
Since the poles of $R(x)$ are given by the set of all the zeros of $M(qx)$ and poles of $r(x)$, the complete description of the poles of (\ref{SolEqInhomogeneousMerCstarNonzero}) follows from Corollaries \ref{HomogeneousMerC} and \ref{HomogeneousMerCstar} and from Proposition \ref{AdditiveDecompositonMCstar}.
\end{proof}

\section{Linear $q$-difference equations of order $n$}\label{ordern}

In this section, we are interested in linear $q$-difference equations of higher order. Using the classical method of factorization \cite[Lemme 9]{zhang}, and, for more details \cite[Section 3.1]{mazh} (see also \cite{BG41,Sau04}), it is well-known that their study can be reduced to that of first order equations. For the sake of readability of the present work, we briefly recall in the next section this method and we illustrate it with two examples.

\subsection{A constructive method for meromorphic solutions}\label{Section31}

Let us consider a $q$-difference equation of order $n\geq2$ of the form
\begin{equation}\label{e688}
\delta y(x)=0,\quad \delta=m_0(x)+m_1(x)\sigma_q+...+m_n(x)\sigma_q^n,
\end{equation}
where $\sigma_q$ stands for the $q$-difference operator $\sigma_qy(x)=y(qx)$, the coefficients $m_j(x)\in\MerC$ are meromorphic functions on $\C$ for all $j=0,...,n$, and where $m_0m_n\not\equiv 0$.

For every $0\le j\le n$, let us set $\Delta_j=\{(j,\hbox{val}(m_j)+t):t\ge0\}\subseteq\R^2$, where $\hbox{val}(m)$ stands for the valuation of $m$ at $0$. Then, defining the Newton polygon associated to (\ref{e688}) as the convex hull of $\bigcup_{0\le j\le n} \Delta_j$, and denoting by $-\infty<k_1\le k_2\le ...\le  k_n<+\infty$ the increasing sequence of its slopes counted with their respective multiplicities (recall that the multiplicity of a slope is the length of the projection of the corresponding edge on the horizontal axis), one can prove that there exists a factorization of the operator $\delta$ of the form
$$m'_0(x)(x^{k_1}\sigma_q-\alpha_1)m'_1(x)(x^{k_2}\sigma_q-\alpha_2)m'_2(x)\cdots (x^{k_n}\sigma_q-\alpha_n)m'_n(x),$$
where $m'_j(x)\in\MerC$ and $\alpha_j\in\C^{\star}$ for all $j=1,...,n$.

The algorithmic procedure for the factorization is described in detail in \cite[Section 3.1]{mazh} (see also \cite{BG41,Sau04}). We illustrate below this method with two examples.

\begin{example}
\textnormal{A first simple example is the linear $q$-difference equation
$$y(q^2x)+(-q^k-1)y(qx)+q^ky(x)=0,$$
for some positive integer $k$. Following the previous algorithm of factorization one has that $n=2$, with $m'_2(x)$ being a constant, $m'_1(x)\equiv m'_0(x)\equiv 1$ and $k_1=k_2=0$, due to the Newton polygon has no slopes. It is straight to check that the factorization of the previous equation is then given by
$$\left(\sigma_q-q^k\right)\left(\sigma_q-1\right)y=0.$$  
We first consider the first order $q$-difference equation $L_1z=0$, with $L_1=\sigma_q-q^k$. From Proposition~\ref{MeromorphicSolutionBasicEquations1} one has that $z(x)=\frac{\Theta_{q}(x)}{\Theta_{q}(q^kx)}=q^{-\frac{k(k+1)}{2}}x^{-k}$ is a solution of this equation. Second, we take equation $L_2y=z(x)$, with $L_2=\sigma_q-1$. The procedure stated in Section~\ref{SectionInhomogeneousLinearQDiffOrder1} determines that $R_{0,\rho}(x)\equiv 0$, $\alpha_{\rho}=0$ and $R_{\infty,\rho}(x)=q^{\frac{k(k+1)}{2}}x^k$. Therefore, we have that $z_{\infty,\rho}(q/x)=\frac{-1}{(q^k-1)q^{\frac{k^2-k}{2}}}\frac{1}{x^k}$. The study of poles and zeros is straight in this example.}
\end{example}

\begin{example}
\textnormal{Let us now consider the second order linear $q$-difference equation
\begin{equation}\label{ExSecondOrderFactTrivial}
y(q^2x)-(qx\sin(2qx)+\cos(x))y(qx)+x\cos(x)\sin(2x)y(x)=0.
\end{equation}
Since it can be factorized into
$$(\sigma_q-\cos(x))(\sigma_q-x\sin(2x))y(x)=0,$$
a meromorphic solution of Eq. (\ref{ExSecondOrderFactTrivial}) is given by a meromorphic solution of the first order inhomogeneous linear $q$-difference equation
$$y(qx)=x\sin(2x)y(x)+r(x),$$
with $r(x)$ a meromorphic solution of $y(qx)=\cos(x)y(x)$.}

\textnormal{Let us now observe that, accordingly to Theorem \ref{HomogeneousHolC}, we can first choose for $r(x)$ an entire function (we have indeed $\cos(x)\in\HolC$ and $\cos(0)=1$). On the other hand, a brief study of the zeros of $x\sin(2x)$ shows that
\begin{itemize}
\item the origin $x=0$ is a double zero of $x\sin(2x)$ with $\frac{x\sin(2x)}{x^2}_{|x=0}=2$;
\item the nonzero zeros of $x\sin(2x)$ are located at the points $n\pi/2$ for all $n\in\Z^*$.
\end{itemize}
Denoting then by $\mathcal P$ the set
$$\mathcal P=\bigcup_{n\in\Z^*}\dfrac{n\pi}{2}q^{\Na},$$
we derive from Corollary \ref{InhomogeneousMerCstar} that, for any $\rho>0$, Eq. (\ref{ExSecondOrderFactTrivial}) admits a meromorphic solution on $\C^*$ whose the poles are located at the elements of
\begin{itemize}
\item$\pm q^{\Z}$ and $2q^{\Z}$:
\item$aq^{\Na^*}$ for any $a\in\mathcal P_{>\rho}$;
\item$aq^{-\Na}$ for any $a\in\mathcal P_{\leq\rho}$.
\end{itemize}
In particular, choosing $\rho\in]0,\pi/2[$, the poles of such a solution are located at the elements of
\begin{itemize}
\item$\pm q^{\Z}$ and $2q^{\Z}$:
\item$\frac{n\pi}{2}q^{\Na^*}$ for any $n\in\Z^*$,
\end{itemize}
and there is no half-spiral of poles with $0$ as accumulation point.}
\end{example}

\subsection{Conclusion and directions for further research}

In Section \ref{sec21}, we have built a meromorphic solution for any first order linear $q$-difference equation with meromorphic coefficients in $\C^*$. Combining then this result with the factorization of linear $q$-difference operator, we have shown in the previous Section \ref{Section31} that this procedure allows to make explicit a meromorphic solution of any linear $q$-difference equation of order $n\geq2$ with meromorphic coefficients in $\C$. However, C. Praagman proved in \cite{pr} that such an equation admits a basis of meromorphic solutions in $\C^*$. So, a possible direction of our further researches is to provide an explicit construction of such a basis.

Another direction of research is related to the factorization of linear $q$-difference operator: in the procedure detailed by F. Marotte and C. Zhang in \cite[Section 3.1]{mazh} to prove the existence of such a factorization, the authors use the fact that the equation under consideration admits an analytic solution, which provides thus a non-constructive proof of the existence of the factorization. Consequently, we can ask the following question: can we explain a constructive algorithm for determining the factorization of any linear $q$-difference operator?


\end{document}